\documentclass[12pt]{article}
\usepackage[centertags]{amsmath}
\usepackage{amsfonts,amssymb,amsthm}
\usepackage{newlfont}

\title{Logic of the Ontological Argument}
\author{{Slobodan Vujo\v{s}evi\'{c}} \\{\small Mathematical Institute}\\[-2mm] {\small SASA,   Belgrade} \and {Filip D. Jevti\'{c}} \\{\small Mathematical Institute}\\[-2mm] {\small SASA,   Belgrade}}

\begin{document}

\maketitle
\begin{abstract}
In his ontological argument G\"{o}del says nothing about its underlying logic. The argument is modal and at least of second-order and since S5 axiom is used so it is widely accepted that the logic of the argument is the S5 second-order modal logic. However, there is a step in the proof in which G\"{o}del applies the necessitation rule on the assumptions of the argument (see \cite{gedel}). This is repeated by all of his followers (see \cite{types} and \cite{hayek}). This application of the necessitation rule can seriously harm the consequence relation of the logic of the ontological argument. It seems that the only way to preserve the modal logic S5 for the ontological argument is to assume some of its axioms in the necessitated form.    
\end{abstract}

\section{Ontological argument}

Ontological argument belongs to the family of arguments that establishes the existence of God by relying only on pure logic. Argument proceeds from the idea of God to the reality of God and was first clearly formulated by St. Anselm in his Proslogion (1077–78). Later famous versions were given by Descartes, Leibniz  and others. These ontological proofs were clearly not formal, in the sense of formal logic, but they show a striking similarity in their form: they argue that God exists (actually, really, necessarily) if God is possible (consistent, present in our mind), and then proceed to prove that God is indeed possible. This form is clearly logical but Kant use the term  "ontological" having in mind ontological context of such proofs. Their formal and modal substance are recognized by  Hartshorne (see \cite{sobel}), who expressed them in S5 modal propositional logic, while G\"{o}del presented his version of formal ontological argument in S5 second-order modal logic. 

G\"{o}del's reasons for his interest in the onthological argument are most clearly expressed in the following quote (see \cite{hayek}): ``I believe that there is much more reason in religion, though not in the churches, than one
commonly believes, but we were brought up from early youth to a prejudgment against
it. We are, of course, far from being able to confirm scientifically the theological world
picture, but it might, I believe, already be possible to perceive by pure reason (without
appealing to the faith in any religion) that the theological worldview is thoroughly
compatible with all known data (including the conditions that prevail on our earth).
The famous philosopher and mathematician Leibniz already tried to do this 250 years
ago, and this is also what I tried in my previous letters (ontological argument).”

\section{Consequence and proof in modal logic}

There is no doubt that the propositional skeleton of the logic of G\"{o}del's argument is the modal logic S5, or something close to it. Besides the axioms of classical propositional logic, the modal axioms of the logic S5 are
\begin{align*}
\Box (A \rightarrow B) &\rightarrow (\Box A \rightarrow \Box B),\\
\Box A &\rightarrow A, \\
\Box A &\rightarrow \Box \Box A \\
\Diamond \Box A &\rightarrow \Box A,
\end{align*}
where $\Box$ is the necessity operator and, $\Diamond$ is the possibility operator defined by $\Diamond A \leftrightarrow \neg \Box \neg A,$ and the inference rules are modus ponens and necessitation: from $A$ infer $\Box A$. The last axiom is usually called the S5 axiom. 

The notions of consequence and proof in modal logic are different from those in classical logic. The relation that a sentence $A$ is a consequence of assumptions $\Sigma$ can have two meanings in modal logic: $A$ is true at each world at which the members of $\Sigma$ are true, and $A$ is true in every model in which $\Sigma$ holds. The two notions are not equivalent and to distinguish between them some authors (see \cite{modal}) are using the terms {\em local} assumption and {\em local} consequence in the first, and the terms {\em global} assumption and {\em global} consequence in the second case. 
We shall show how this semantical distinction is reflected in the syntax of modal logic.

Assume a sentence $A$ globally; if $A$ is true at an arbitrary world $w$ in some model, then $A$ is true at every world accessible to $w$ (since $A$ holds in every world), so $\Box A$ holds at $w$. Since $w$ is arbitrary, $\Box A$ holds at every world of a model. This means that if $A$ is a global assumption, the necessitation rule can be applied to $A$.
On the other hand, if we assume $A$ locally, so that $A$ is known to be true at a world $w$ of some model, there is no reason to expect that $\Box A$ is also true at $w$. If $A$ is a local assumption, the necessitation rule cannot apply to it. 

The distinction between global and local assumptions in formal deductions comes down to the applicability or nonapplicability of the necessitation rule. A formal proof or derivation in modal logic does not allow the use of the necessitation rule to local premises and their consequences. To insure this, some authors define modal derivations as finite sequences divided in two separate parts, global and local (see \cite{modal}). The global part comes first, containing only global premises and the necessitation rule is allowed, while the local part comes second containing local premises, but without the necessitation rule.

\section{Necessitation in G\"{o}del's argument}

It is well known that G\"{o}del was involved in the foundation of the modern approach to modal logic. He was among the first logicians who introduced the necesitation rule that made possible the simple and elegant modal axiom systems that are in use today. But in the early 1970s, at the time G\"{o}del wrote his note about the ontological argument, the idea of possible world semantics was new and perhaps not well appreciated. G\"{o}del argument is modal and is presented in at least second-order logic, however the exact logic is not specified.

According to what we have told about consequence in modal logic, to allow the unrestricted use of the necessitation rule in the logic S5, we have to assume the axioms of our theories globally. In the modal logic S5, where $\Box A \leftrightarrow \Box \Box A$, assuming $A$ globally we assume $\Box A$. Formally, this means that all axioms of the theories in the S5 logic must come in the necessitated form, i.e. with $\Box$ prefixed.

G\"{o}del's argument is a particular version of the general ontological argument that usually means two things: to prove that God's existence is possible and to prove that God exists necessarily if it exists. If $Q$ is the statement that God exists, this means that in the general ontological argument we have to prove $\Diamond Q$ and $Q \rightarrow \Box Q$ (Anselm's principle). It is generally accepted that with these assumptions within S5 logic one can prove $\Box Q$: the necessitation of Anselm's principle gives $\Diamond Q \rightarrow \Diamond \Box Q$, the S5 axiom gives $\Diamond Q \rightarrow \Box Q$ and the first assumption finally gives $\Box Q$ (see \cite{types}, \cite{gedel}, \cite{hayek} and \cite{sobel}). But the use of necessitation in this proof was not correct. It seems that the only way to overcome this incorrectness is to formulate Anselm's principle in the form $\Box (Q \rightarrow \Box Q)$: it is necessary that God exists necessary if He exists. 

At some point in his note, relying on axioms that are not formulated in necessitated form, G\"{o}del presents the theorem $$G(x) \rightarrow \Box \exists y G(y), $$ where $G(x)$ means that ``$x$ is godlike being" (see page 403 in \cite{gedel}), and without any comments proceedes in the following three steps:
\begin{align*}
\exists x G(x) &\rightarrow \Box \exists y G(y), \\ 
\Diamond \exists x G(x) &\rightarrow \Diamond \Box \exists y G(y), \\
\Diamond \exists x G(x) &\rightarrow \Box \exists y G(y).
\end{align*}
In the first step the existential quantifier is introduced, the second step comes from the necessitation rule, and the third uses of the S5 axiom. Since he was able to prove $\Diamond  \exists x G(x)$, G\"{o}del finally concludes $\Box \exists y G(y)$. 

G\"{o}del, as well as his followers and commentators in this matter, say nothing about the local or global character of the ontological argument axioms. They present these axioms in the unnecessitated form (see \cite{modal}, \cite{hayek} and \cite{sobel}), and use the necessitation rule on them and on their consequences. Perhaps they have in mind global axioms? 

\bigskip

\noindent {\em Acknowledgement.} We express our thanks to Kosta Do\v sen who warned us about the problem concerning the correctness of the necessitation rule in the ontological argument.


\begin{thebibliography}{1}

\bibitem{types}  Fitting, M., 2002, {\em Types, Tableaux and G\"{o}del’s God},  Dordrecht, Kluwer.

\bibitem{modal} Fitting, M. and Mendelsohn, R. L., 1990, {\em First-Order Modal Logic}, Dordrecht, Kluwer.

\bibitem{gedel} G\"{o}del, K., 1995, {\em Ontological Proof}, In {\em Collected Works, Volume III, Unpublished essays and lectures}, New York, Oxford University Press, pp. 403-404. 

\bibitem{hayek} H\' ayek, P. 2011,   {\em G\"{o}del's Ontological Proof and Its Variants}, In {\em Kurt G\"{o}del and the Foundations
of Mathematics  - Horizons of Truth},  Cambridge, Cambridge University Press, 2011, pp. 307-324.

\bibitem{sobel} Sobel, J. H., 1987,  {\em G\"{o}del’s Ontological Proof}, In {\em On Being and Saying: Essays in Honor of Richard Cartwright},  Cambridge, MA, MIT Press, pp. 24-61.

\end{thebibliography}
\end{document}